\documentclass[preprint,12pt]{elsarticle}

\usepackage{amssymb}
\usepackage{stmaryrd}
\usepackage{amsfonts,amsmath,latexsym}
\usepackage{mathrsfs}
\usepackage{amsbsy}
\usepackage{indentfirst}
\usepackage{amsmath}
\usepackage{amsfonts}
\usepackage{bbm}
\usepackage{float} 
\usepackage{mathrsfs}
\usepackage{dsfont}
\usepackage{graphicx}
\usepackage{subfigure}%
\usepackage{rotating}
 \def\esssup{\mathop {\rm ess\,sup}}
  \def\essinf{\mathop {\rm ess\,inf}}

\numberwithin{equation}{section}

\newtheorem{theorem}{Theorem}[section]

\newtheorem{remark}{Remark}[section]
\newtheorem{definition}{Definition}[section]

\newcommand{\beqnar}{\begin{eqnarray*}}
\newcommand{\eeqnar}{\end{eqnarray*}}
\newcommand{\ba}{\begin{array}}
\newcommand{\ea}{\end{array}}

\journal{}

\allowdisplaybreaks[4]

\begin{document}

\begin{frontmatter}



\title{From classical to modern central limit theorems }

 \author{Vladimir V. Ulyanov
 }
 
  \address{Lomonosov Moscow State University and HSE University
  }
  


  \begin{abstract}
De Moivre (1733), investigating the limit distribution of the binomial distribution, was the first to discover the existence of the normal distribution and the central limit theorem. In this review article, we briefly recall the history of classical central limit theorem and martingale central limit theorem, and introduce a new direction of central limit theorem, namely nonlinear central limit theorem and nonlinear normal distribution.
   \end{abstract}


\end{frontmatter}



\section{The Classical Central Limit Theorems}
 The central limit theorem (CLT) is one of the gems of probability theory. Its importance can hardly be overestimated both from the theoretical point of view and from the point of view of applications in various fields. 
The first version of the CLT was appeared as the De Moivre–Laplace theorem.  De Moivre's  investigation was motivated by a need to compute probabilities of winning in various games of chance. In the proof De Moivre used Stirling's formula, which he knew.
\begin{theorem} (De Moivre (1733))
 Let $\{X_n\}_{n\ge 1}$ be a sequence of independent Bernoulli random variables, each with success probability $p\in(0,1)$, that is, for each $i$,
$$P(X_i=1)=p=1-P(X_i=0).$$
Let $S_n=\sum_{i=1}^{n}X_i$ denote the  cumulative number of successes in first $n$ Bernoulli trials.  Then for any $a<b\in\mathbb{R}$,
\begin{equation}
	P\left(a< \frac{S_n-np}{\sqrt{np(1-p)}}\leq b\right)\to \Phi(b)-\Phi(a), \nonumber
\end{equation}
where $$\Phi(x)=\frac{1}{\sqrt{2\pi}}\int_{-\infty}^xe^{- {y^2}/{2}}dy$$ is the cumulative distribution function of the standard normal distribution.	
\end{theorem}
Thus De Moivre discovered the probability distribution, which in the late 19th century came to be called the normal distribution. Another name -- Gaussian distribution -- is used in honor of Gauss, who arrived at this distribution as ‘‘the law of error’’ in his famous work Gauss(1809) on problems of measurement in astronomy and the least squares method.

Laplace proved the Moivre-Laplace theorem anew using the Euler- McLaurin summation formula.
 \begin{theorem} (Laplace (1781)) In the notation of Theorem 1.1, for any $a\in\mathbb{R}$, we have
\begin{eqnarray}
   P\left(|{S_n-np - z}|\leq a\right) = 2 \left(\Phi\left(\frac{a\sqrt{n}}{\sqrt{x\,x'}}\right)-\Phi(0)\right) + \frac{\sqrt{n}}{\sqrt{2\pi x'x}}\exp\left(-\frac{a^2n}{2x'x}\right),\nonumber
\end{eqnarray}
     where  $z\in\mathbb{R}$, $|z|<1$ and $x=np+z,\,\, x'=n(1-p)-z$.
 \end{theorem}
This theorem, besides convergence to the normal law, also gives a good estimate of the accuracy of the normal approximation.

 
After the first CLT, many famous mathematicians studied the CLT, such as Possion,  Dirichlet, Cauchy,  and so on.  
The first significant generalization of the Moivre-Laplace theorem was the Lyapunov theorem:
\begin{theorem} (Lyapunov(1900, 1901))  Let $\{X_n\}_{n\ge 1}$ be a sequence of independent  random variables 
	with mean $E[X_n]=\mu_n$, variance  $Var[X_n]=\sigma_n^2$ and finite moment $E[|X_n - E X_n|^{2+\delta}], \,\, \delta > 0. $. Let $S_n=\sum_{i=1}^{n}X_i$, $B_n=\sum_{i=1}^{n}\sigma_i^2$. If   for some $\delta>0$, 
	\begin{equation}
		\label{lc}
		\lim_{n\to \infty}\frac{1}{B^{1+\delta/2}_n}\sum_{i=1}^{n}E[|X_i-E[X_i]|^{2+\delta}  =0. \nonumber
	\end{equation}
 Then for any $a<b\in\mathbb{R}$,
 \begin{equation}
 	P\left(a< \frac{S_n-E[S_n]}{\sqrt{B_n}}\leq b\right)\to \Phi(b)-\Phi(a), \nonumber
 \end{equation}
uniformly with respect to $a$ and $b$.
\end{theorem}
This theorem was proved by a new method: the method of characteristic functions. The formulation of the problem and a possible solution were proposed by Chebyshev (1887), who suggested to use the method of moments  by comparing the moments of sums of independent random variables with the moments of the Gaussian distribution. The   method of moments is still helpful in some   cases nowadays.  


  In the twenties of the last century, the study of CLT introduced modern probability theory. The most crucial CLT in the early years of the   century belongs to L\'evy.  In 1922, L\'evy's fundamental theorems on characteristic functions were proved.

  At the same time, Lindeberg also used characteristic functions to study CLT, and a fundamental CLT, called the Lindeberg-Feller CLT, was formulated:

\begin{theorem} In the notation of Theorem 1.3, for 
	\begin{equation}
		\label{fc}
		\lim_{n\to \infty}\frac{\max_{1\le i\le n}\sigma_i^2}{B_n}=0 \nonumber
	\end{equation}
and for any $b$ 
\begin{equation}
	P\left(\frac{S_n-E[S_n]}{\sqrt{B_n}}\leq b\right)\to \Phi(b),\,\,\,\,n\to\infty, \nonumber
\end{equation}
to hold, it is necessary and sufficient that the condition (the Lindeberg condition) for any $\varepsilon>0$:
   \begin{equation}
\label{lc}
\lim_{n\to \infty}\frac{1}{B_n}\sum_{i=1}^{n}E[(X_i-E[X_i])^2 1_{\{|X_i-E[X_i]|>\varepsilon  \sqrt{B_n}\}}]=0 \nonumber
\end{equation}
is met.
%
\end{theorem} 
Sufficiency was proved by Lindeberg(1922) , necessity by Feller(1935).

%


A detailed and extensive review of classical CLT can be found in Fischer(2011).
Further research was carried out in several directions, among them: obtaining estimates of the accuracy of the normal approximation, see e.g. Petrov(1995), G\"{o}tze et al.(2017), Shevtsova (2014); generalization of CLT to the multivariate case, see e.g. Bhattacharya and Ranga Rao (1976), Sazonov(1981), and infinite dimensional case, see e.g. Bentkus et al.(1991), G\"{o}tze and Zaitsev (2014), Prokhorov and Ulyanov(2013).

    The classical CLT mainly studies the sums of independent random variables.  It is pretty different to derive the limiting distribution of the sums of dependent random variables.

   \section{The Martingale Central Limit Theorems}

  The classical result is that independent and identically distributed variables lead to a normal distribution under the proper moment condition. Counterexamples show that violation of the independence or identity of the distribution leads to a non-normal limit distribution. However, numerous examples also show that violation of independence or identity of distribution can still lead to a normal distribution. 
   Bernstein (1926) and L\'evy (1935) independently put forward a new direction of study: how to prove the CLT for sums of dependent random variables. In 1935, L\'evy established a CLT under some conditions, which can be regarded as the beginning of the martingale CLT.

    In the classical CLT,  $S_n=\sum_{i=1}^nX_i$,  $\{X_n\}_{n\ge 1}$ is assumed as a sequence of independent random variables. In the contents of the martingale CLT, $S_n=\sum_{i=1}^nX_i$ is a martingale,  $\{X_n\}_{n\ge 1}$ is assumed as a sequence of martingale difference sequences, i.e.   $E[X_n|\mathcal{F}_{n-1}]=0$,  where $\{\mathcal{F}_{n}\}_{n\ge 0}$ is a sequence of given $\sigma$-algebras filtration.  In 1935,  L\'evy established his martingale CLT, which reads as follows.

 \begin{theorem} (L\'evy (1935))
     Let $\{X_n\}_{n\ge 1}$ be a sequence of  random variables defined on $(\Omega, \mathcal{F}, \{\mathcal{F}_{n}\}_{n\ge 0}, P)$, with   $E[X_n|\mathcal{F}_{n-1}]=0$. Denote $S_n=\sum_{i=1}^nX_i$,  $\sigma_n^2=E[X_n^2|\mathcal{F}_{n-1}]$ and $b_n^2=\sum_{i=1}^{n}\sigma_{i}^2$.  If
   \begin{equation}
\label{l1}
\sum_{n=1}^{\infty}\sigma_n^2=\infty, ~~P\text{-a.s.},  \nonumber
\end{equation}
and for any $\varepsilon>0$,
 \begin{align}
&\lim_{n\to \infty}\sum_{i=1}^{n}P(|X_i|>\varepsilon b_n|\mathcal{F}_{i-1})=0, &P\text{-a.s.},
  \nonumber\\
&\lim_{n\to \infty}\frac{1}{b_n}\sum_{i=1}^{n}E\left[X_i1_{\{|X_i|>\varepsilon b_n\}}|\mathcal{F}_{i-1}\right]=0, &P\text{-a.s.},
\nonumber\\
&\lim_{n\to \infty}\frac{1}{b_n^2}\sum_{i=1}^{n}E\left[X_i^21_{\{|X_i|>\varepsilon b_n\}}|\mathcal{F}_{i-1}\right]=0, &P\text{-a.s.},
\nonumber\\
&\lim_{n\to \infty}\frac{1}{b_n^2}\sum_{i=1}^{n}\left(E\left[X_i1_{\{|X_i|>\varepsilon b_n\}}|\mathcal{F}_{i-1}\right]\right)^2=0, &P\text{-a.s.}.
\nonumber
 \end{align}
    Then
                  \begin{equation}
\label{clt2}
\frac{S_n}{\sqrt{n}}\stackrel{d}{\longrightarrow}  N(0,c^2), \nonumber
\end{equation}
where $c$ is a constant.
 \end{theorem}

L\'evy's result depends on the condition for $\sigma_n^2$, which are random variables.  The assumptions of L\'evy's theorem are too strict.  Many authors tried to relax the assumptions of
   L\'evy's result, including Doob (1953), Billingsley (1961), Ibragimov (1963), Cs\"org\"o (1968),  and so on.  These works led to the  CLT under some conditions, including
       \begin{equation}
\label{l6}
\frac{b_n^2}{E[S_n^2]}\stackrel{P}{\longrightarrow}  C, \nonumber
\end{equation}
  where $C$ is a constant.

  Brown (1971) improved the previous martingale CLT and obtained the following result.

\begin{theorem} (Brown (1971))
    Let $\{X_n\}_{n\ge 1}$ be a sequence of  random variables defined on $(\Omega, \mathcal{F}, \{\mathcal{F}_{n}\}_{n\ge 0}, P)$, with  $E[X_n|\mathcal{F}_{n-1}]=0$. Denote $S_n=\sum_{i=1}^nX_i$,  $\sigma_n^2=E[X_n^2|\mathcal{F}_{n-1}]$,  $\varphi_n(t)=E[e^{itX_n}|\mathcal{F}_{n-1}]$ $b_n^2=\sum_{i=1}^{n}\sigma_{i}^2$, and $s_n^2=E[S_n^2]$.  If
  \begin{equation}
\label{l7}
\frac{b_n^2}{s_n^2}\stackrel{P}{\longrightarrow}  1, \nonumber
\end{equation}
 \begin{equation}
\label{l8}
\prod_{j=1}^{n}\varphi_j\left(\frac{t}{\sqrt{s_n^2}}\right)\stackrel{P}{\longrightarrow}  e^{-\frac{t^2}{2}},
\end{equation}
and
 \begin{equation}
\label{l9}
\frac{\max_{1\le i \le n}\sigma_i^2}{s_n^2}\stackrel{P}{\longrightarrow}  0.    \nonumber
\end{equation}
 Then
                  \begin{equation}
\label{clt2-1}
\frac{S_n}{\sqrt{s_n^2}}\stackrel{d}{\longrightarrow}  N(0,1).   \nonumber
\end{equation}

\end{theorem} 
The condition (\ref{l8}) is a critical condition.  McLeish (1974) introduced an elegant method with a new martingale CLT.

\begin{theorem} (McLeish (1974))
Let $\{X_{n,i}, n\ge 1, 1\le i \le k_n\} $ be an array of  random variables defined on $(\Omega, \mathcal{F}, P)$,  with $E_P[X_{n,i}|\mathcal{F}_{n,i-1}]=0$, where $\mathcal{F}_{n,i}=\sigma(X_{n,j}, 1\le j \le i)$. Denote $S_n=\sum_{i=1}^{k_n}X_{n,i}$, $T_n=\prod_{i=1}^{k_n}(1+itX_{n,i})$.  If for all  real $t$,  $\{T_n\}$ is uniformly integrable, and
    \begin{equation}
\label{l19}
\lim_{n\to \infty}E[T_n]=1,~~~\sum_{i=1}^{k_n}X_{n,i}^2\stackrel{P}{\longrightarrow}  1, \nonumber \end{equation}
 \begin{equation}
\label{l10}
\max_{1\le i \le k_n} |X_{n,i}|\stackrel{P}{\longrightarrow}  0.  \nonumber
\end{equation}
 Then
                  \begin{equation}
\label{clt3}
S_n \stackrel{d}{\longrightarrow}  N(0,1).  \nonumber
\end{equation}
\end{theorem} 

  McLeish's elegant method proves the condition (\ref{l8}) holds.  The condition (\ref{l8}) implies that the limiting distribution of martingales $S_n$ is Gaussian.   However,  condition (\ref{l8})  means that the conditional characteristic function converges to a constant. It is an unnatural condition.  The limiting distribution will differ when the conditional characteristic function or conditional variance converges to a random variable.  Hall (1977) obtained the following result, an essential progress in the martingale CLT.

 \begin{theorem} (Hall (1977))
    Let $\{X_{n,i}, n\ge 1, 1\le i \le k_n\} $ be an array of  random variables defined on $(\Omega, \mathcal{F}, P)$,  with $E[X_{n,i}|\mathcal{F}_{n,i-1}]=0$, where  $\mathcal{F}_{n,i}=\sigma(X_{j,k_n}, 1\le j \le i)$ and for $n\ge 1$, $1\le i\le k_n$,  $\mathcal{F}_{n,i} \subset \mathcal{F}_{n+1,i}$, Denote $S_n=\sum_{i=1}^{k_n}X_{n,i}$.  If     \begin{equation}
\label{l11}
\lim_{n\to \infty}E[\max_{1\le i \le k_n} X_{n,i}^2]=0,  \nonumber \end{equation}
 \begin{equation}
\label{l12}
\sum_{i=1}^{k_n} X_{n,i}^2\stackrel{P}{\longrightarrow}  T.
\end{equation}
 Then
                  \begin{equation}
\label{clt4}
S_n \stackrel{d}{\longrightarrow}  T'N(0,1),             \nonumber
\end{equation}
  where $T'$ is a independent copy of $T$. 
 \end{theorem} 

  Hall's result implies that the limiting distribution of martingale may not be a Gaussian. It may be a conditional Gaussian.  In the study of economics and finance,  condition (\ref{l12})  usually fits the real-life situation. More often than not, final decisions can only be made in an ambiguous context. We are often faced with decision-making problems under uncertainty in real-life situations. In most cases,  the classical CLT and normal distribution are not suitable.  For example, we can not directly construct the confidence interval utilizing the martingale CLT.

  The theory of nonlinear probabilities and expectations has developed rapidly over the last thirty years and became an important tool in investigating model uncertainty or ambiguity. The nonlinear CLT is a significant research area that describes the asymptotic behavior of a sequence of random variables with distribution uncertainty. Its limiting distribution is no longer represented by the classical normal distribution but by a series of "new nonlinear normal distributions", such as the $g$-expectation and $G$-normal distribution.
  The nonlinear CLT in nonlinear probability and expectation theory can fill the vast gaps between the martingale CLT and real life.

\section{Nonlinear Central Limit Theorems}
There are two main frameworks for studying nonlinear CLT. The nonlinear expectation framework $(\Omega,\mathcal{H},\mathbb{E})$ proposed by Peng is one approach to characterize distributional uncertainty. Another approach involves a set of probability measures $\mathcal{P}$ on $(\Omega,\mathcal{F})$ to study the nonlinear CLT, as used by Chen, Epstein and their co-authors.

\subsection{Nonlinear CLT under nonlinear expectations}
  In 1997,  Peng constructed a large class of dynamically consistent nonlinear expectations through backward stochastic differential equations, known as $g$-expectation, with the corresponding dynamic risk measure referred to as $g$-risk measure. The $g$-expectation can handle uncertain probability sets  $\left\{P_\theta, \theta \in \Theta\right\}$ controlled by a given probability measure $P$. However, for singular cases (i.e. $P(A)=0$ while $P_\theta(A)>0$), $g$-expectation is no longer applicable. 
  Peng, breaking free from the framework of the original probability space, created the theory of nonlinear expectation spaces
    and introduced a more general nonlinear expectation $G$-expectation, see Peng (2019).

\begin{definition} Given a set $\Omega$, let $\mathcal{H}$ be a linear space of real-valued functions defined on $\Omega$. If the functional $\mathbb{E}$: $\mathcal{H} \rightarrow  \mathbb{R}$ satisfies the following four conditions:

\noindent 1), Monotonicity: If $X \geqslant Y$, then $\mathbb{E}[X] \geqslant \mathbb{E}[Y]$;\\
2), Preserving Constants: $\mathbb{E}[c]=c$, for all $c \in \mathbb{R}$;\\
3), Subadditivity: $\mathbb{E}[X+Y] \leqslant \mathbb{E}[X]+\mathbb{E}[Y]$, for all $X, Y \in \mathcal{H}$;\\
4), Positive Homogeneity: $\mathbb{E}[\lambda X]=\lambda \mathbb{E}[X]$, for all $\lambda \geqslant 0$.

Then, the functional $\mathbb{E}$ is called a {\bf sublinear expectation}. The triple $(\Omega, \mathcal{H}, \mathbb{E})$ is referred to as a {\bf sublinear expectation space}. If only conditions 1) and 2) are satisfied, $\mathbb{E}$ is termed a {\bf nonlinear expectation}, and $(\Omega, \mathcal{H}, \mathbb{E})$ is called a {\bf nonlinear expectation space}. 
\end{definition}

In the framework of nonlinear expectations, starting from fundamental assumptions, concepts such as the distribution of random variables, independence, correlation, stationarity, Markov processes, and so forth, can also be derived. Simultaneously, nonlinear Brownian motion and the corresponding stochastic analysis represent substantial extensions of classical stochastic analysis theory. Furthermore, limit theorems still hold under nonlinear expectations.  Peng (2008) developed an elegant partial differential equation method and obtained the first nonlinear CLT under sublinear expectations.

\begin{theorem} (Peng (2008))
    Let $\left\{X_i\right\}_{i\ge1}$ be a sequence of independent and identically distributed  random variables in the sublinear expectation space $(\Omega, \mathcal{H}, \mathbb{E})$. Assume that
$$
\mathbb{E}\left[X_1\right]=\mathbb{E}\left[-X_1\right]=0, \quad \lim _{c \rightarrow \infty} \mathbb{E}\left[\left(\left|X_1\right|^2-c\right)^{+}\right]=0 .
$$
Let $S_n:=\sum_{i=1}^nX_i$. Then, for any $\varphi\in C(\mathbb{R})$ with linear growth condition,
\begin{equation}
    \lim _{n \rightarrow \infty} \mathbb{E}\left[\varphi\left(\frac{S_n}{\sqrt{n}}\right)\right]=\mathbb{E}[\varphi(\xi)],
\end{equation}
where $\xi\sim \mathcal{N}(0,[\underline{\sigma}^2,\overline{\sigma}^2])$ is a $G$-normally distributed random variable and the corresponding sublinear function $G:  \mathbb{R} \mapsto \mathbb{R}$  is defined by 
$$
G(a):=\mathbb{E}\left[ \frac{a}{2}X_1^2 \right],   a\in \mathbb{R} .
$$
\end{theorem}   

\subsection{Nonlinear CLT under a set of probability measures}
Peng's nonlinear CLT opens up a new field to replace the martingale CLT in the study of economics and finance.    In economic markets, random variables objectively exist, but the uncertain probability measure $P$ measurement may not be deterministic. The well-known Ellsberg paradox illustrates that random variables exist, but finding a probability measure $P$ to quantify a given random variable is not always possible. This example demonstrates that in practical applications of probability theory, there may be ambiguity in people's understanding of probability measures, and it needs to be clarified which measure should be used to quantify uncertainty better. The economic community often refers to this as ambiguity,  indicating the uncertainty arising from the market and people's limited cognitive abilities. From the perspective of the development of economic markets, economists have found that Kolmogorov's established probability axioms can quantify the intrinsic laws of economic market development. However, it cannot precisely characterize the external impact of human behavior on market development laws. Thus, the classical and martingale CLT does not work. Therefore, a new discipline called behavioral economics has emerged to study the laws of economic markets using different perspectives (probability measures). Establishing universal asymptotic results is a significant and challenging issue in this situation.

Inspired by this, Chen, Epstein and their co-authors, unlike Peng who started with nonlinear expectations, investigated nonlinear CLT  for a set of probability measures (i.e., in the context of ambiguity). They use a set of probability measures $\mathcal{P}$, which is assumed be "rectangular" (or closed with respect to the pasting of alien marginals and conditionals), to describe the distribution uncertainty  of the random variables $\{X_n\}$ defined on $(\Omega,\mathcal{F})$. Given the history information $\{\mathcal{G}_{n-1}\}$, the conditional mean and variance of $\{X_n\}$  will vary under different probability measures $Q\in\mathcal{P}$, which leads to the concepts of upper and lower conditional means and variances, see (\ref{mubar}) and (\ref{variance-up-low}).
Thus, they focused on two characteristics, mean uncertainty and variance uncertainty, respectively, to get a nonlinear CLT. 

Particularly, Chen, Epstein and their co-authors found and established two types of
nonlinear normal distributions, which have the explicit probability densities as given in (\ref{proba-density}) and (\ref{tranproba}), to characterize the 
 limiting distribution in the nonlinear CLT. These explicit expressions are the first explicit formulas for nonlinear CLT since de Moivre (1733), Laplace (1781), and Gauss (1809) discovered and proved the (linear) normal distribution for a single probability measure more than two hundred years ago. It is worth noting that these two types of non-linear normal distributions also play a crucial role in studying multi-armed bandits and quantum computing as well as  nonlinear statistics ( see e.g. Chen et al.(2023)).  
\medskip

\noindent{\bf Case I: CLT with mean uncertainty}
\medskip

 Chen and Epstein (2022) established a family of CLT with mean uncertainty under a set of probability measures. Under the assumptions of constant conditional variance of the random variable sequence and conditional mean constrained to a fixed interval $[\underline{\mu},\bar{ \mu}]$, they proved that the limiting distribution can be described by the $g$-expectation or a solution of a backward stochastic differential equation (BSDE). Furthermore, for a class of symmetric test functions, they show the limiting distribution has an explicit density function, given in (\ref{proba-density}).


\begin{theorem} (Chen and Epstein (2022))   Let  $(\Omega, \mathcal{F})$ be a measurable space, $\mathcal{P}$ be a family of probability measures on  $(\Omega, \mathcal{F})$, and $\left\{X_i\right\}$ be a sequence of real-valued random variables defined on this space. 
The history information is
represented by the filtration $\{\mathcal{G}_{i}\}_{i\ge1}$, $({\mathcal{G}_{0}=\{\emptyset,\Omega\}})$, such that $\left\{  X_{i}\right\}  $ is adapted to
$\{\mathcal{G}_{i}\}$.

Assume that the upper and lower conditional means of   $\{X_{i}\}$ satisfy:
\begin{equation}
\esssup\limits_{Q\in\mathcal{P}}E_{Q}[X_{i}|\mathcal{G}_{i-1}]=\overline{\mu}\text{ and }\essinf_{Q\in\mathcal{P}}E_{Q}[X_{i}|\mathcal{G}_{i-1}%
]=\underline{\mu}\text{, for all }i\geq1\text{.} \label{mubar}%
\end{equation}
Assume that $\{X_{i}\}$ has an \emph{unambiguous conditional variance
}$\sigma^{2}$, that is,
\begin{equation}
E_{Q}\left[  (X_{i}-E_{Q}[X_{i}|\mathcal{G}_{i-1}])^{2}|\mathcal{G}%
_{i-1}\right]  =\sigma^{2}>0\text{ for all }Q\in\mathcal{P}\text{ and all
}i\text{.} \label{condvar}%
\end{equation} 
Further more, assume that $\mathcal{P}$ is rectangular and $\left\{X_i\right\}$ satisfies the Lindeberg condition:
\begin{equation}\label{linder}
    \lim _{n \rightarrow \infty} \frac{1}{n} \sum_{i=1}^n \sup_{Q \in \mathcal{P}} E_Q\left[\left|X_i\right|^2 I_{\left\{\left|X_i\right|>\sqrt{n \varepsilon}\right\}}\right]=0, \forall \varepsilon>0 .
\end{equation}
Then for all $\varphi \in C([-\infty, \infty])$,
\small
\begin{equation}
\lim\limits_{n\rightarrow\infty}\sup_{Q\in\mathcal{P}}E_{Q}\left[
\varphi\left(  \frac{1}{n}{\sum_{i=1}^{n}X_{i}}+\frac{1}{\sqrt{n}}%
\sum\limits_{i=1}^{n}\frac{1}{\sigma}{(X_{i}-E_{Q}[X_{i}|\mathcal{G}_{i-1}%
])}\right)  \right]  =\mathbb{E}_{\left[  \underline{\mu},\overline{\mu
}\right]  }[\varphi\left(  B_{1}\right)  ],\label{CLT}%
\end{equation}
or equivalently,%
\begin{equation}
\lim\limits_{n\rightarrow\infty}\inf_{Q\in\mathcal{P}}E_{Q}\left[
\varphi\left(  \frac{1}{n}{\sum_{i=1}^{n}X_{i}}+\frac{1}{\sqrt{n}}{\sum
_{i=1}^{n}\frac{1}{\sigma}(X_{i}-E_{Q}[X_{i}|\mathcal{G}_{i-1}])}\right)
\right]  =\mathcal{E}_{\left[  \underline{\mu},\overline{\mu}\right]
}[\varphi\left(  B_{1}\right)  ],\label{CLTinf}%
\end{equation}\normalsize
where $\mathbb{E}_{\left[  \underline{\mu},\overline{\mu}\right]  }%
[\varphi\left(  B_{1}\right)  ]\equiv Y_{0}$ \textbf{is called $g$-expectation
by Peng (1997),} given that $(Y_{t},Z_{t})$ is the solution of the BSDE%
\begin{equation}
Y_{t}=\varphi(B_{1})+\int_{t}^{1}\max\limits_{\underline{\mu}\leq\mu
\leq\overline{\mu}}(\mu Z_{s})ds-\int_{t}^{1}Z_{s}dB_{s},\ 0\leq t\leq1, \nonumber
\label{BSDE-1}%
\end{equation}
and $\mathcal{E}_{\left[  \underline{\mu},\overline{\mu}\right]  }%
[\varphi\left(  B_{1}\right)  ]\equiv y_{0}$, given that $(y_{t},z_{t})$ is
the solution of the BSDE
\begin{equation}
y_{t}=\varphi(B_{1})+\int_{t}^{1}\min\limits_{\underline{\mu}\leq\mu
\leq\overline{\mu}}(\mu z_{s})ds-\int_{t}^{1}z_{s}dB_{s},\ 0\leq
t\leq1\text{.}   \nonumber\label{BSDE-2}%
\end{equation}
Here $(B_{t})$ is a standard Brownian motion

Particularly, when $\varphi$ is symmetric with center $c\in\mathbb{R}$, that is, $\varphi(c+x)=\varphi(c-x)$, and is monotonic on $(c,\infty)$, the limits of (\ref{CLT}) and (\ref{CLTinf}) can be expressed explicitly.  

\noindent {\bf (1)} If $\varphi$ is increasing on $(c,\infty)$, then 
\begin{equation} \footnotesize \lim\limits_{n\rightarrow\infty}\sup_{Q\in\mathcal{P}}E_{Q}\left[
\varphi\left(  \frac{1}{n}{\sum_{i=1}^{n}X_{i}}+\frac{1}{\sqrt{n}}%
\sum\limits_{i=1}^{n}\frac{1}{\sigma}{(X_{i}-E_{Q}[X_{i}|\mathcal{G}_{i-1}%
])}\right)  \right]=\int_{\mathbb{R}} \varphi(y)f^{\frac{\overline\mu-\underline\mu}{2},\frac{\overline\mu+\underline\mu}{2},c}(y)dy\label{clt-mean-explicit1}
\end{equation}
\begin{equation} \footnotesize \lim\limits_{n\rightarrow\infty}\inf_{Q\in\mathcal{P}}E_{Q}\left[
\varphi\left(  \frac{1}{n}{\sum_{i=1}^{n}X_{i}}+\frac{1}{\sqrt{n}}%
\sum\limits_{i=1}^{n}\frac{1}{\sigma}{(X_{i}-E_{Q}[X_{i}|\mathcal{G}_{i-1}%
])}\right)  \right]=\int_{\mathbb{R}} \varphi(y)f^{\frac{\underline\mu-\overline\mu}{2},\frac{\overline\mu+\underline\mu}{2},c}(y)dy\label{clt-mean-explicit2}
\end{equation}
\noindent {\bf (2)} If $\varphi$ is decreasing on $(c,\infty)$, then 
\begin{equation}  \footnotesize \lim\limits_{n\rightarrow\infty}\sup_{Q\in\mathcal{P}}E_{Q}\left[
\varphi\left(  \frac{1}{n}{\sum_{i=1}^{n}X_{i}}+\frac{1}{\sqrt{n}}%
\sum\limits_{i=1}^{n}\frac{1}{\sigma}{(X_{i}-E_{Q}[X_{i}|\mathcal{G}_{i-1}%
])}\right)  \right]=\int_{\mathbb{R}} \varphi(y)f^{\frac{\underline\mu-\overline\mu}{2},\frac{\overline\mu+\underline\mu}{2},c}(y)dy\label{clt-mean-explicit3}
\end{equation}
\begin{equation} \footnotesize \lim\limits_{n\rightarrow\infty}\inf_{Q\in\mathcal{P}}E_{Q}\left[
\varphi\left(  \frac{1}{n}{\sum_{i=1}^{n}X_{i}}+\frac{1}{\sqrt{n}}%
\sum\limits_{i=1}^{n}\frac{1}{\sigma}{(X_{i}-E_{Q}[X_{i}|\mathcal{G}_{i-1}%
])}\right)  \right]=\int_{\mathbb{R}} \varphi(y)f^{\frac{\overline\mu-\underline\mu}{2},\frac{\overline\mu+\underline\mu}{2},c}(y)dy\label{clt-mean-explicit4}
\end{equation}
where the density function $f^{\alpha,\beta,c}$ is given as follows:
\begin{equation}
f^{\alpha,\beta,c}(y)= \frac{1}{\sqrt{2\pi }}e^{-\frac{(y-\beta)^2-2\alpha (|y-c|-|c-\beta|)+\alpha^2}{2}}-\alpha e^{2\alpha |y-c|}
\Phi(-|c-\beta|-|y-c|-\alpha),\label{proba-density}
\end{equation}
The above density function of the Chen-Epstein distribution degenerates into the density function of the classical normal (Gaussian) distribution only when $\alpha=0$.
\end{theorem} 
\begin{remark} Let  $\beta=0$ and $c=0$, the density function $f^{\alpha,\beta,c}$ has the following properties:
\begin{itemize}
 \item If $\alpha <0$, the curve of $f^{\alpha,\beta,c}$ is more spike than the normal distribution, referred as a spike distribution. When we use $f^{\alpha,\beta,c}$ to denote maximum probability density in (\ref{clt-mean-explicit3}) and the minimum probability density in (\ref{clt-mean-explicit2}), the corresponding $\alpha<0$.
 
 \item If $\alpha>0$, the curve of $f^{\alpha,\beta,c}$ is similar to two normal distributions  hand in hand, referred as a binormal distribution. When we use $f^{\alpha,\beta,c}$ to denote maximum probability density in (\ref{clt-mean-explicit1}) and the minimum probability density in (\ref{clt-mean-explicit4}), the corresponding $\alpha>0$.
 \item If $\alpha=0$, the curve of $f^{\alpha,\beta,c}$ is degenerated to a standard normal distribution. From (\ref{clt-mean-explicit1}) to (\ref{clt-mean-explicit4}), when $\overline{\mu}=\underline\mu$ the correspionding $\alpha=0$.
\end{itemize}
The density function of a Bandit distribution is shown in the following figures.
\end{remark}

The curves of the density function $f^{\alpha,\beta,c}$, for the cases that $\beta=0,c=0$ and  $\alpha$ taken different values, are shown in the following figures.
 \begin{figure}[H]
  \centering
    \includegraphics[width=4.5in]{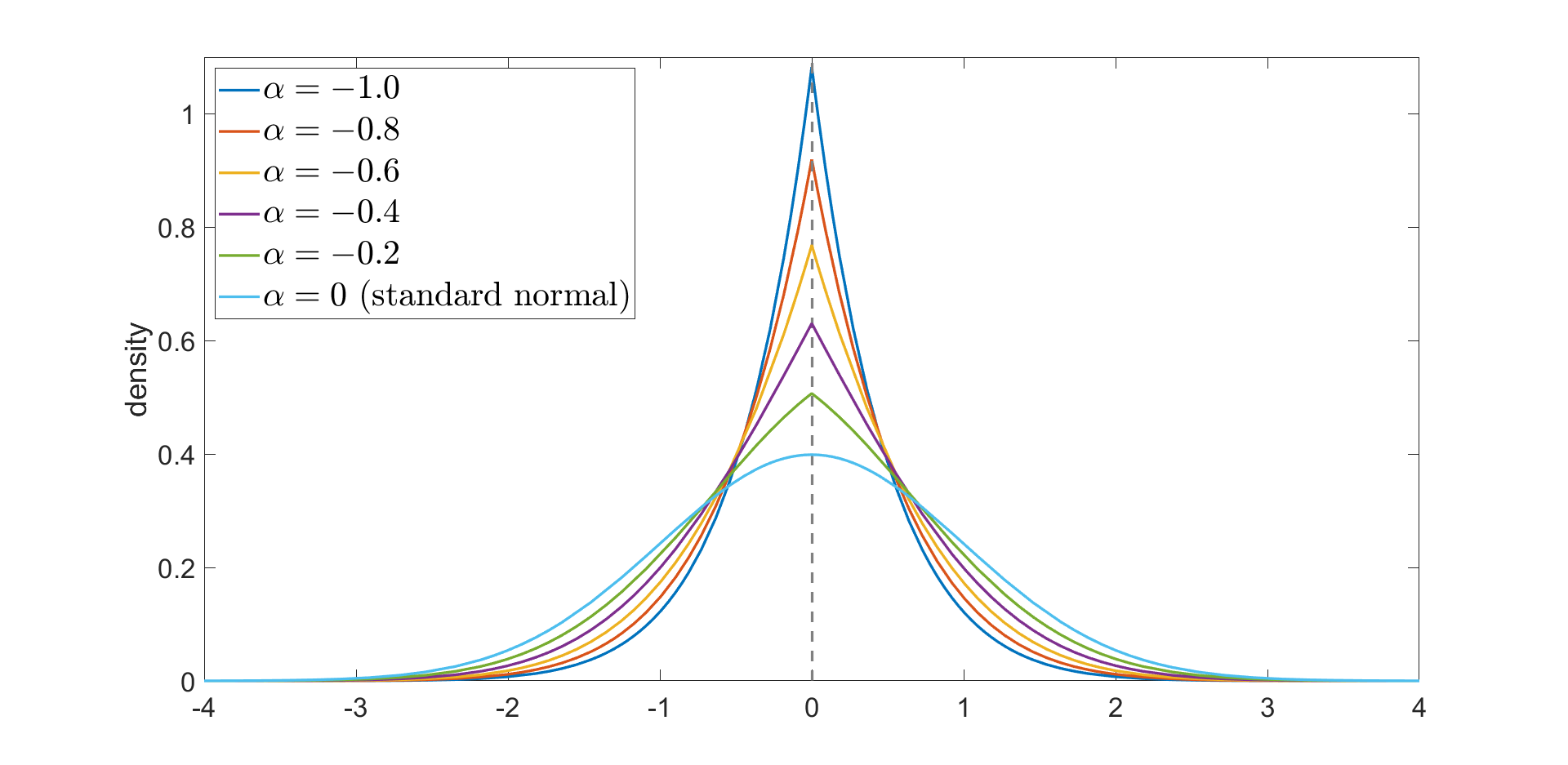}
    \caption{$\beta=0,c=0$ and $\alpha\le0$}
\end{figure}
 \begin{figure}[H]
  \centering
    \includegraphics[width=4.5in]{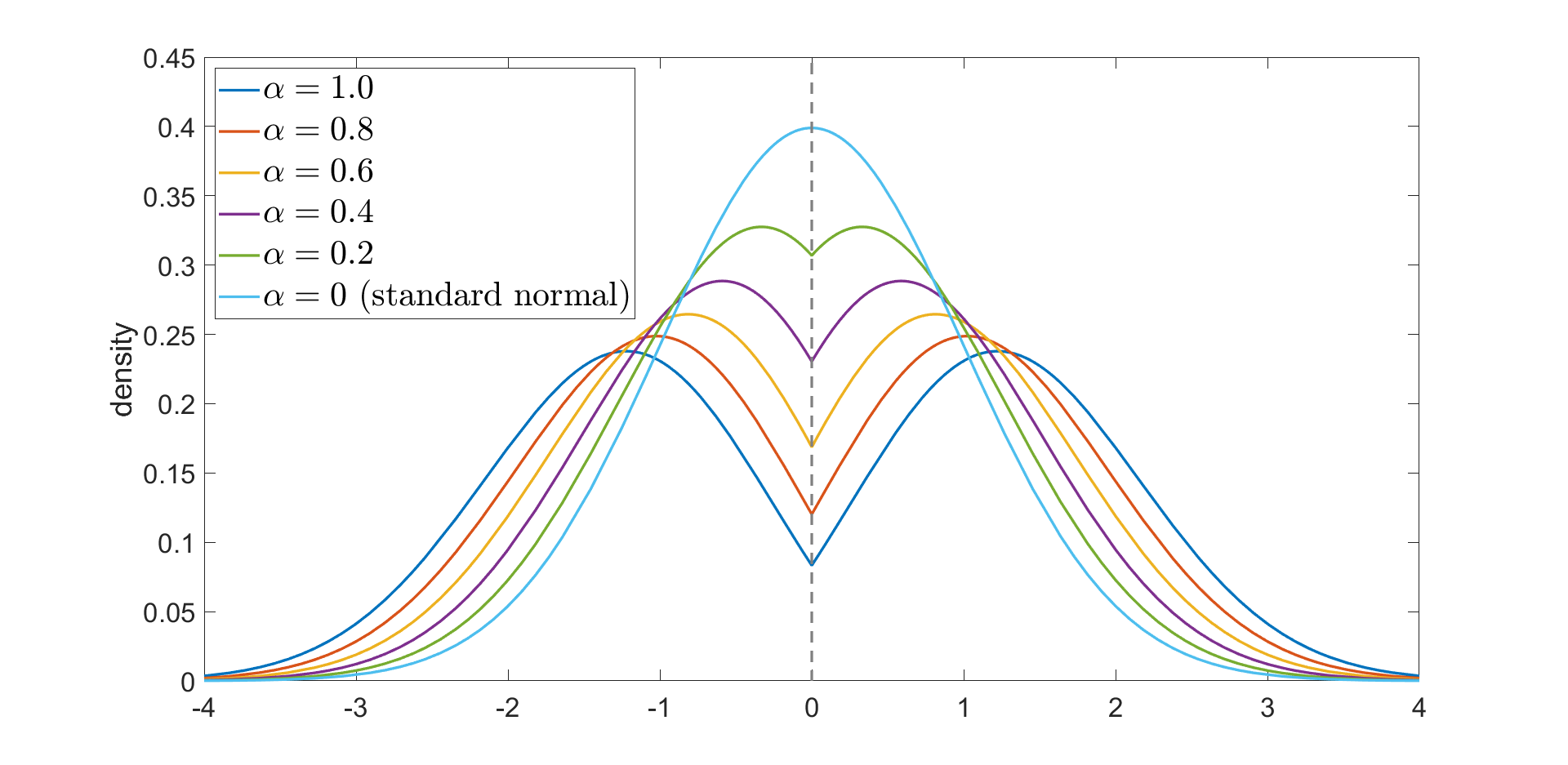}
    \caption{$\beta=0,c=0$ and $\alpha\ge0$}
\end{figure}

\medskip

\noindent{\bf Case II: CLT with variance uncertainty}
\medskip

Chen, Epstein and Zhang (2023) investigated the CLT with variance uncertainty  under a set of probability measures. The considered random variables sequence has an unambiguous conditional mean and its conditional variance is constrained to vary within the interval $\left[\underline{\sigma}^2, \bar{\sigma}^2\right]$. First of all, the limiting distirbution can still be described by $G$-normal distribution. More significantly, for a class of "S-Shaped" test functions, which are important indexes for characterizing loss aversion in behavioral economics, they demonstrated that the limiting distribution also has an explicit probability density function, given in (\ref{tranproba}).

\begin{theorem}\label{nclt-variance} ( Chen, Epstein and Zhang (2023))
Let  $(\Omega, \mathcal{F})$ be a measurable space, $\mathcal{P}$ be a family of probability measures on  $(\Omega, \mathcal{F})$, and $\left\{X_i\right\}$ be a sequence of real-valued random variables defined on this space. The history information is
represented by the filtration $\{\mathcal{G}_{i}\}_{i\ge1}$, $({\mathcal{G}_{0}=\{\emptyset,\Omega\}})$, such that $\left\{  X_{i}\right\}  $ is adapted to
$\{\mathcal{G}_{i}\}$.

Assume that $\{X_{i}\}$ has an \emph{unambiguous conditional mean
}$0$, that is,
\begin{equation}
E_{Q}\left[  X_{i}|\mathcal{G}%
_{i-1}\right]  =\mu=0\text{ for all }Q\in\mathcal{P}\text{ and all
}i\text{,}                      \nonumber\label{condmean}%
\end{equation} 
Assume that the upper and lower conditional variances of $\{X_{i}\}$ satisfy:
\begin{equation}
\esssup\limits_{Q\in\mathcal{P}}E_{Q}[X_{i}^2|\mathcal{G}_{i-1}]=\overline{\sigma}^2\text{ and }ess\inf_{Q\in\mathcal{P}}E_{Q}[X_{i}^2|\mathcal{G}_{i-1}%
]=\underline{\sigma}^2\text{, for all }i\geq1\text{.} \label{variance-up-low}%
\end{equation}
Assume also the Lindeberg condition (\ref{linder}) is met 
and that $\mathcal{P}$ is rectangular.  Put $S_{n}%
=\sum_{i=1}^{n}X_{i}$. For any $\varphi\in C([-\infty,\infty])$, we have
\begin{equation}\label{CLT-G}
\lim_{n\to\infty}\sup_{Q\in\mathcal{P}}E_{Q}\left[\varphi\left(\frac{S_n}{\sqrt{n}}\right)\right]=\mathbb{E}[\varphi(\xi)]  \nonumber
\end{equation}
where $\xi\sim \mathcal{N}(0,[\underline\sigma^2,\overline\sigma^2])$ is a $G$-Normal distribution under the sublinear expectation space $(\Omega,\mathcal{H},\mathbb{E})$.

Particularly, for
any $c\in\mathbb{R}$ and $\varphi_{1}\in C_{b}^{3}\left(  \mathbb{R}\right)
$, set $\theta=\underline{\sigma}/\overline{\sigma}$ and define functions
\begin{equation}
\phi(x)=\left\{
\begin{array}
[c]{ll}%
\varphi_{1}(x) & x\geq c\\
\varphi_{2}(x)=-\theta\varphi_{1}\left(  -\frac{1}{\theta}(x-c)+c\right)
+\left(  1+\theta\right)  \varphi_{1}(c) & x<c  
\end{array}
\right. \nonumber \label{function-ph-theta}%
\end{equation}%
\begin{equation}
\overline{\phi}(x)=\left\{
\begin{array}
[c]{ll}%
\varphi_{1}(x) & x\geq c\\
\overline{\varphi}_{2}(x)=-\frac{1}{\theta}\varphi_{1}\left(  -\theta
(x-c)+c\right)  +\left(  1+\frac{1}{\theta}\right)  \varphi_{1}(c) & x<c
\end{array}
\right.                        \nonumber\label{function-phbar}%
\end{equation}

\begin{description}
\item[(1)] If $\varphi_{1}^{\prime\prime}(x)\leq0$ for $x\geq c$,
then%
\begin{align}
\lim_{n\rightarrow\infty}\sup_{Q\in\mathcal{P}}E_{Q}\left[  \overline{\phi}\left(
\frac{S_{n}}{\sqrt{n}}\right)  \right]   &  
=\int_{\mathbb{R}}\overline{\phi}(y)q^{\underline{\sigma},\overline{\sigma},c}(y)dy,\label{thm-clt-eq-O-leq-sup}\\
\lim_{n\rightarrow\infty}\inf_{Q\in\mathcal{P}}E_{Q}\left[  \phi\left(  \frac{S_{n}}%
{\sqrt{n}}\right)  \right]   &  
=\int_{\mathbb{R}}\phi(y)q^{\overline{\sigma},\underline{\sigma},c}(y)dy.\label{thm-clt-eq-O-leq-inf}%
\end{align}

\item[(2)] If $\varphi_{1}^{\prime\prime}(x)\geq0$ for $x\geq c$,
then%
\begin{align}
\lim\limits_{n\rightarrow\infty}\sup_{Q\in\mathcal{P}}E_{Q}\left[  \phi\left(  \frac{S_{n}%
}{\sqrt{n}}\right)  \right]   &  
=\int_{\mathbb{R}}\phi(y)q^{\overline{\sigma},\underline{\sigma},c}(y)dy,\label{thm-clt-eq-O-geq-sup}\\
\lim_{n\rightarrow\infty}\inf_{Q\in\mathcal{P}}E_{Q}\left[  \overline{\phi}\left(
\frac{S_{n}}{\sqrt{n}}\right)  \right]   & 
=\int_{\mathbb{R}}\overline{\phi}(y)q^{\underline{\sigma},\overline{\sigma},c}(y)dy.\label{thm-clt-eq-O-geq-inf}%
\end{align}

 
\end{description}

The density function of $q^{\alpha,\beta,c}$ is given as follows:
\begin{align}
q^{\alpha,\beta,c}(y)=   \frac{1}{\sqrt{2\pi}\sigma
(y)}e^{-\frac{\left(  \frac{c}{\sigma(0)}+\frac{y-c}{\sigma
(y)}\right)  ^{2}}{2}}    +\tfrac{\beta-\alpha}{\beta+\alpha}\tfrac
{sgn(y-c)}{\sqrt{2\pi}\sigma(y)}e^{-\frac{\left(  \left|  \frac{c}{\sigma
(0)}\right|  +\left|  \frac{y-c}{\sigma(y)}\right|  \right)  ^{2}}%
{2}}  ,\quad\label{tranproba}%
\end{align}
and  $\sigma(y)=\alpha I_{[c,\infty)}(y)+\beta I_{(-\infty,c)}(y),\ \forall
y\in\mathbb{R}$.

The above density function degenerates into the density function of the  normal distribution $N(0,\sigma^2)$ only when $\alpha=\beta=\sigma$.
\end{theorem}
\begin{remark}
    The assumptions and results of Theorem \ref{nclt-variance} is slightly different  from those in the formal publication of Chen, Epstein and Zhang (2023). 
    The frameworks and assumptions in the formal publication are slightly more complicated  due to the consideration of the bandit problem and the construction of the optimal strategies in the Bayesian framework. 
\end{remark}

The following figure shows the curves of  the density  $q^{\alpha,\beta,c}$.
 \begin{figure}[H]
  \centering
    \includegraphics[width=4.5in]{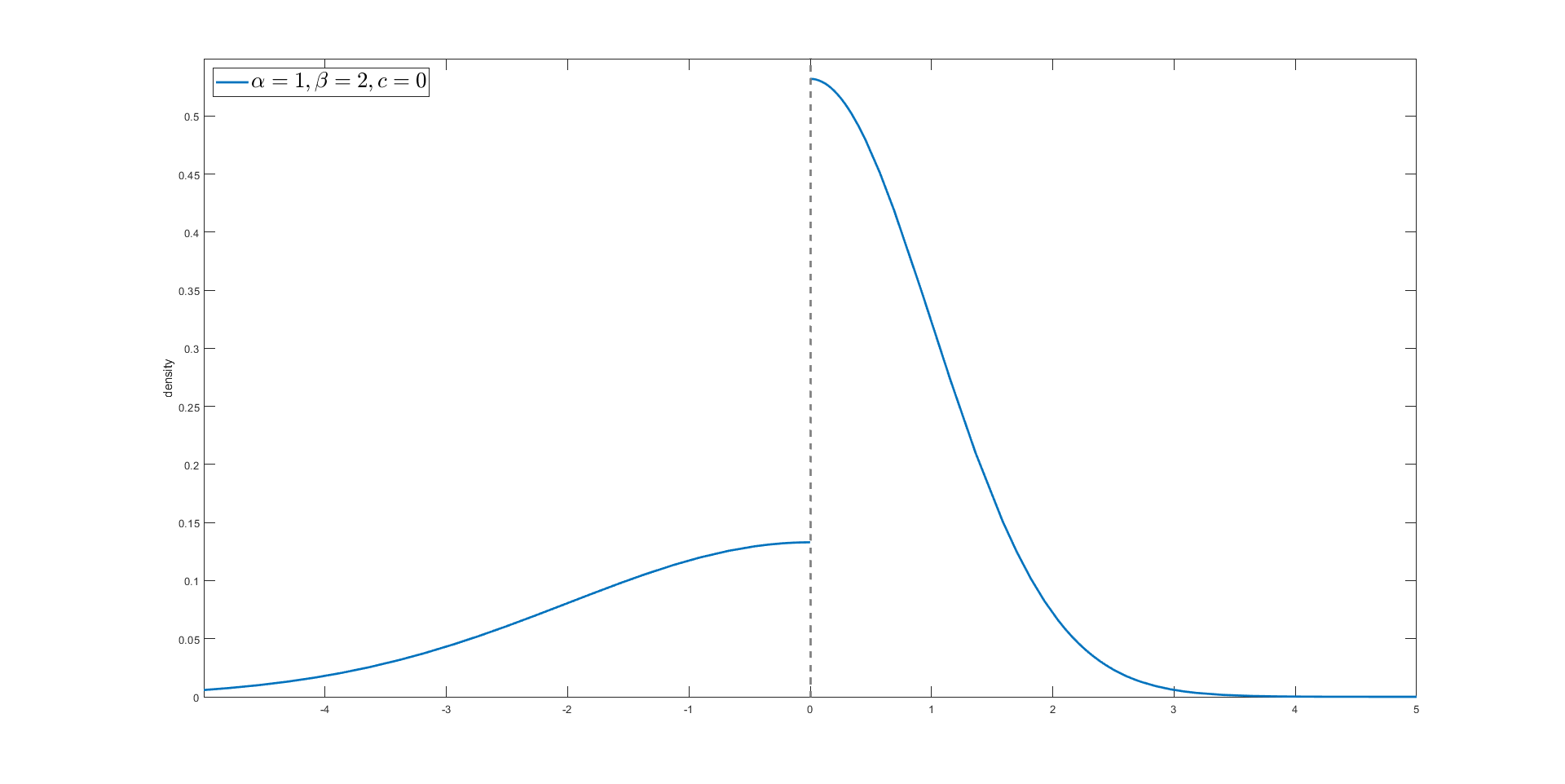}
    \caption{The maximum probability density of (\ref{thm-clt-eq-O-leq-sup}) or the minimum probability density of (\ref{thm-clt-eq-O-geq-inf})}
    \includegraphics[width=4.5in]{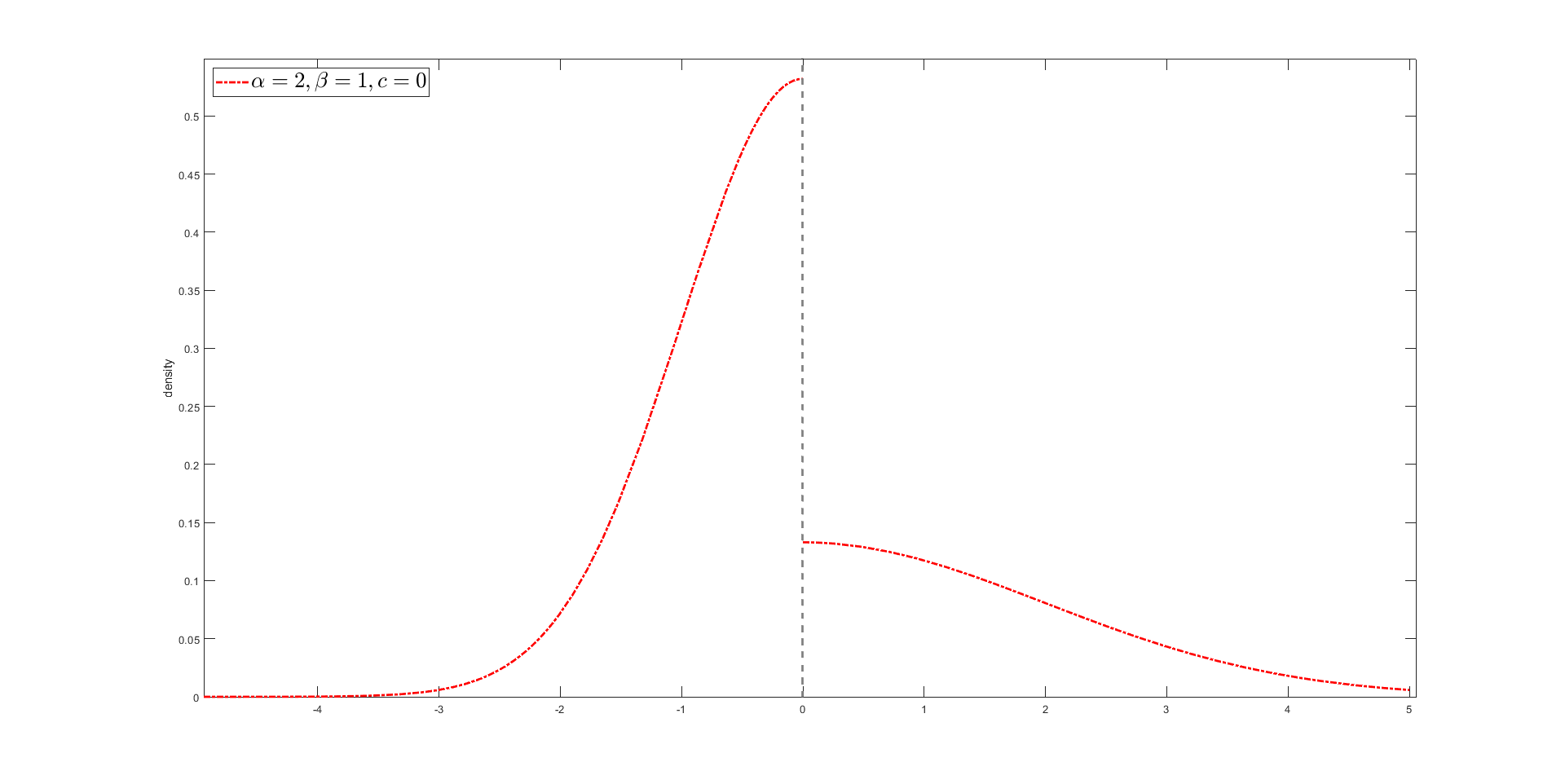}
    \caption{The minimum probability density of (\ref{thm-clt-eq-O-leq-inf}) or the maximum probability density of (\ref{thm-clt-eq-O-geq-sup})}
\end{figure}
 When compare the curves of classical normal distribution, one has
 \begin{figure}[H]
  \centering
    \includegraphics[width=5.5in]{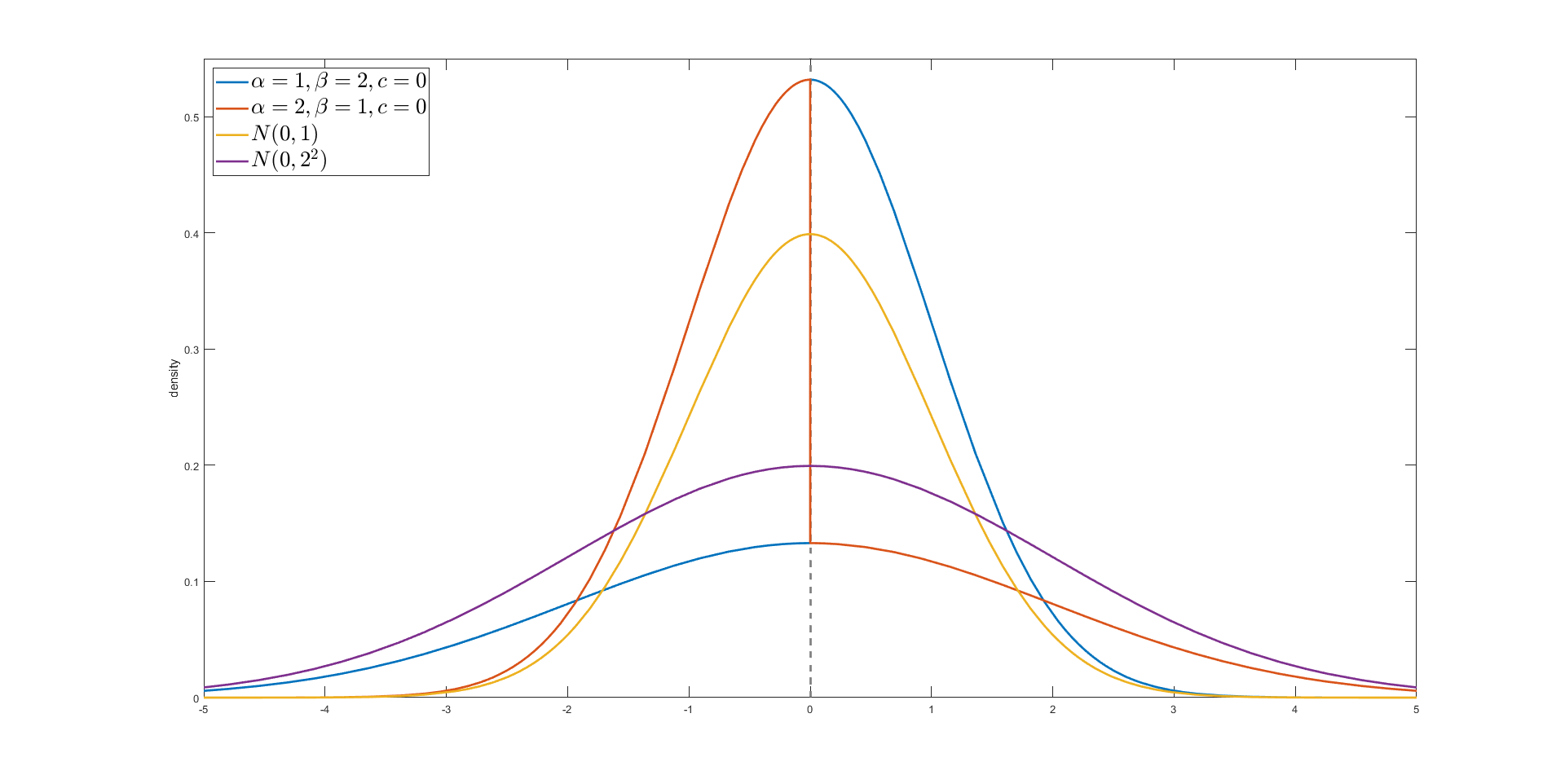}
\end{figure}

\section{Differences between  classical CLT and Nonlinear CLT}
For convenience, we use \textbf{CLT} to represent  the classical CLT , \textbf{NE-CLT} represent the nonlinear Peng's CLT on nonlinear expectations given by Peng, and  \textbf{NP-CLT} represent the nonlinear Chen-Epstein CLT under a set of probability measures given by Chen and Epstein as well as Zhang.
\begin{sidewaystable}\footnotesize
     \centering
    \begin{tabular}{|c|c|c|c|}\hline
         & CLT & NE-CLT & NP-CLT\\ \hline
        Frameworks & $(\Omega,\mathcal{F},P) $ & $ (\Omega,\mathcal{H},\mathbb{E}) $ & $
             (\Omega,\mathcal{F},\mathcal{P})$\\ \hline
        Independence & $\begin{array}{c}
             \{X_i\} \text{ is independent} \\
             \text{ or martingale difference} 
        \end{array}$ & $\begin{array}{c}
             \{X_i\} \text{ is independent} \\
             \text{under }(\Omega,\mathcal{H},\mathbb{E})\text{ i.e.}\\
             \mathbb{E}[f(X_1,\cdots,X_n)]\\
             =\mathbb{E}\left[\mathbb{E}[f(x,X_n)]_{x=(X_1,\cdots,X_{n-1})}\right] 
        \end{array}$ & $\begin{array}{c}
              \text{No concept of independence,} \\
             \mathcal{P}\text{ is rectangular} \\
\sup_{Q\in\mathcal{P}}E_Q\left[f(X_1,\cdots,X_n)\right]\\
=\sup_{Q\in\mathcal{P}}E_Q\left[\esssup_{Q\in\mathcal{P}}E_Q\left[f(X_1,\cdots,X_n)|\mathcal{G}_{n-1}\right]\right]
        \end{array}$ \\ \hline
        $\begin{array}{c}
              \text{Mean and} \\
            \text{Variance} 
        \end{array}$  & $\begin{array}{c}
              \text{Identical distributed} \\
            \text{Same mean and variance} 
        \end{array}$  & $\begin{array}{c}
              \text{Upper and lower means} \\
              \overline{\mu}=\mathbb{E}[X_1],\underline{\mu}=-\mathbb{E}[-X_1]\\
            \text{Upper and lower variances} \\
               \overline{\sigma}^2=\mathbb{E}[X_1^2],\underline{\sigma}^2=-\mathbb{E}[-X_1^2]
        \end{array}$  & $\begin{array}{c}\text{Upper and lower conditional means } (\ref{mubar}) \\
            \text{Upper and lower variances } (\ref{variance-up-low})
        \end{array}$ \\ \hline
      Limit distribution   & Normal distribution & $\begin{array}{c}G\text{-normal distribution } \\
          G\text{-distribution}  
        \end{array}$  & $\begin{array}{c}\text{Solution of BSDE } \\
            G\text{-normal distribution } \\
          \text{Explicit nonlinear normal distribution}  
        \end{array}$ \\ \hline
        Proofs & $\begin{array}{c}\text{Characteristic functions} \\
             \text{Moment method} \\
          \text{Stein’s method}  \\
          \text{The Lindeberg exchange method}
        \end{array}$ & $\begin{array}{c}G\text{-heat equation combine with} \\
             \text{Taylor's expansion} 
        \end{array}$ & $\begin{array}{c}\text{Guess the form of the limiting distribution} \\
             \text{Construct basic functions} \{H_t(x)\} \\
           \text{such as }H_t(x)=\mathbb{E}[\varphi(x+\sqrt{1-t}\xi)]\\
           \text{Prove the nice properties of }H_t\\
           \text{Taylor's expansion with the idea of}\\
        \text{Lindeberg exchange method}
        \end{array}$\\ \hline
    \end{tabular}
    \caption{Differences between  classical CLT and Nonlinear CLT}
    \label{tab:my_label}
\end{sidewaystable}

\end{document}